\documentstyle[12pt]{article}
\begin{document}
\null\vskip1cm
FT-425-1997
\vskip1cm
\centerline{\Large\bf q-Deformed Schr\" odinger Equation}
\vskip1.52cm
\centerline{{\large M. Micu}
\footnote{E-mail address:~MICUM@THEOR1.IFA.RO}}
\centerline{Department of Theoretical Physics}
\centerline{ Horia Hulubei Institute of Physics and Nuclear
Engineering}
\centerline{ POB MG-6, Bucharest, 
76900 Romania}
\vskip2cm
{\bf Abstract} We found hermitian realizations of the position vector
$\vec{r}$, angular momentum $\vec{\Lambda}$ and linear momentum
$\vec{p}$ behaving like vectors with respect to the $SU_q(2)$
algebra, generated by $L_0$ and $L_\pm$. They are used to write the $q$
deformed Schr\" odinger equation, whose solution for Coulomb and
oscillator potential are briefly discussed. 
\newpage
{\bf I. INTRODUCTION}
\vskip0.5cm
The quantum mechanics of a point-like particle is
constructed starting with two vectors: the position
vector, $\vec{r}$, and the linear momentum, 
$\vec{p}=-ih{\partial\over\partial\vec{r}}$, having the well
known commutation relations. 
These two vectors are used to build all the other
quantities, like the angular momentum, the interacting
potential, etc. In general, these new quantities are
noncommutative, their commutation relations being determined by
the commutation relations satisfied by $\vec{r}$ and $\vec{p}$. 

In a $q$-deformed quantum mechanics the
commutation relations between the generators of the $SU_q(2)$
algebra, $\vec{L}$, and the position vector, $\vec{r}$, are
well defined and it is natural to take these vectors as the basic quantities 
from which all the other must be built. 

Wishing to build a $q$-deformed Schr\" odinger Hamiltonian we
searched for a realization of $\vec{p}$ entering the kinetic
energy term. First it was necessary to find a realization
for $\vec{r}$ 
and for $\vec{L}$ as self adjoint quantities obeying the known 
commutation relations. Then one has to look for a realization of  
$\vec{p}$ in terms of $\vec{r}$ and of $\vec{L}$. We found that
$\vec{p}$ can be written as a sum of two terms which are respectively
parallel and perpendicular to $\vec{r}$. The first one
is assumed to have the simplest form and is written as
$-i{\vec{r}\over r^2}~(r{\partial\over\partial r}+1)$
while the second one is expressed as a vector product of $\vec{r}$ and of
$\vec{L}$. 

The paper is organized as follows: Section II contains the
general commutation relations involving the $q$-angular momentum
and some quantities having definite
transformation properties with respect to the $SU_q(2)$ algebra,
like the invariants C, c and the vector $\vec{\Lambda}$. In the
third section we give a realization of the position vector,
$\vec{r}$, and of the $q$-angular momentum$^{1,2,3}$, $\vec{L}$, in terms of
the polar coordinates $r,~x_0={\rm cos}\theta,~\varphi$.
The realization of the linear momentum $\vec{p}$ is given in the
fourth section. We first build the part perpendicular to
$\vec{r}$, denoted $\vec{\partial}$, using the cross product
$\vec{r} \times\vec{\Lambda}$ and find that it satisfies$^1$ some
commutation relations similar to those satisfied by $\vec{r}$.
The part parallel to $\vec{r}$, supposed to have the simplest
form, is just that coming from the ordinary ${\partial\over
\partial\vec{r}}$. Section V contains the eigen functions of the
$q$-angular momentum$^{2,3}$, written like series in $x_0={\rm
cos}\theta$. The result is a generalization of
the hypergeometric functions $_2{\rm F}_1(a,b,c;{1\over2};x_0^2)$
and  $_2{\rm F}_1(a,b,c;{3\over2};x_0^2)$ which can be related to
the $q$-deformed spherical functions $Y_{lm}(q,x_0,\varphi)$.
Some properties and relations satisfied by the eigen functions
are also listed. In the last section the $q$-deformed
Schr\" odinger equation$^{1,4,5}$ with scalar potential is given.
Its solutions for Coulomb and three dimensional oscillator
potentials are briefly discussed. 
\vskip1cm
{\bf II. THE {\large$q$}-ANGULAR MOMENTUM}
\vskip0.5cm
The $SU_q(2)$ algebra is generated by
three operators $L_+,~L_-$ and $L_0$, also named the $q$-angular
momentum, having the following commutation relations:
$$\left[~L_0~,~L_\pm~\right]~=~\pm~L_\pm\eqno(1)$$
$$\left[~L_+~,~L_-~\right]~=~\left[2~L_0\right]\eqno(2)$$
where the quantity in square brackets is defined as
$$\left[n\right]~=~{q^n-q^{-n}\over q-q^{-1}}.\eqno(3)$$
In the following we shall introduce quantities having definite
transformation properties with respect to
the $SU_q(2)$ algebra. They will be further used to build
$q$ scalars and $q$ vectors, like, for instance, the $q$ linear
momentum, entering the expression of the hamiltonian operator.

First of all we remind that $SU_q(2)$ algebra has an invariant,
C, called the Casimir operator 
$${\rm C}~=~L_-~L_+~-~\left[L_0\right]~\left[ L_0~+1\right
]\eqno(4)$$
whose eigenvalue in the $(2l+1)$ dimensional irreducible
representation is:
$${\rm C}_l=[l]~[l+1].\eqno(5)$$

A vector in this algebra is a set of three quantities
$v_k,~k=\pm1, 0$ satisfying the following relations:
$$\left[~L_0~,~v_k~\right]~=~k~v_k\eqno(6)$$
$$\left[~L_\pm~v_k~-~q^k~v_k~L_\pm~\right]~q^{L_0}~=~\sqrt{[2]}~
v_{k\pm1}\eqno(7)$$
where $v_{\pm2}$ must be set equal to zero in the right hand side
for $k=\pm1$.

By comparing the relations (1), (2) with (6), (7), we observe
that, unlike the $SU(2)$ algebra, the 
operators $L_k$ do not represent the components of a vector in
the above sense. 
However, one can use $L_\pm$ and $L_0$ to define a vector,
$\vec{\Lambda}$, in the following manner:
$$\Lambda_{\pm1}~=~{\mp 1\over\sqrt{[2]}}~q^{-L_0}~L_\pm\eqno(8)$$
$$\Lambda_0~=~{1\over[2]}~\left(q~L_+~L_-~-~q^{-1}~L_-~L_+\right).
\eqno(9)$$
It is now an easy matter to show that $\Lambda_k$  satisfy the
relations (6) and (7) as required.

Two vectors $\vec{u}$ and $\vec{v}$ can be used to build a
scalar, $S$, according to the following definition:
$$S~=~\vec{u}~\vec{v}~=~-{1\over q}~u_1~v_{-1}~+~u_0~v_0~-~q~
u_{-1}~v_{1}.\eqno(10)$$
By introducing a generalization of the cross product, the
vectors can be used also to build a new vector, as it will be 
further shown. 

In the case $\vec{u}=\vec{v}=\vec{\Lambda}$, the scalar product
$\vec{\Lambda}^2$ defines a second invariant$^6$, ${\rm C}'$, which
is not independent of C. The eigenvalue of ${\rm C}'$ is
$${\rm C}'_l~=~{[2l]\over[2]}~{[2l+2]\over[2]}\eqno(11)$$
In this paper a third invariant, c, defined as
$${\rm c}~=~q^{-2L_0}~+~\lambda~\Lambda_0\eqno(12)$$
with
$$\lambda~=~q-{1\over q}\eqno(13)$$
will be frequently used in order to write the formulae in a more
compact form.
Its eigenvalue is:
$${\rm c}_l~=~{q^{2l+1}~+~q^{-2l-1}\over[2]}.\eqno(14)$$
It is worth noticing that, in the limit $q=1$ the first and
second invariants ${\rm C,~C}'$ both go into the Casimir
invariant C=$\vec{L}^2=l(l+1)$, while the third one, c,
becomes equal to unity.
The results listed in this section are valid for any realization
of the $SU_q(2)$ algebra. 
\vskip1cm
{\bf III. THE POSITION VECTOR {\large$\vec{r}$} AND A REALIZATION
OF {\large $L_\pm$} AND {\large$L_0$}}
\vskip0.5cm
In $R_q(3)$ space the position vector $\vec{r}$ has three
noncommutative components $r_1,~r_{-1}$ and $r_0$, satisfying the
following relations:
$$r_0~r_{\pm1}~=~q^{\mp2}~r_{\pm1}~r_0\eqno(15)$$
$$r_1~r_{-1}~=~r_{-1}~r_1~+~\lambda~r_0^2.\eqno(16)$$
The quantity $r^2$ defined as
$$r^2~=~\vec{r}~^2~=~-{1\over q}~r_1~r_{-1}~+~r_0^2~-~q~
r_{-1}~r_1\eqno(17)$$
commute with all $r_i$ and with all $L_i$ if
$\vec{r}$ satisfies the conditions (6) and (7) to be a vector. For
$q=1$ the scalar $r$ is nothing else than the length of the
position vector $\vec{r}$. We shall keep this meaning also for
$q\neq1$. 

Searching for the concrete realization of $\vec{r}$,
$L_\pm$ and of $L_0$, we begin by expressing $L_0$ like in 
$R(3)$ case: 
$$L_0~=~-i~{\partial\over\partial\varphi}.\eqno(18)$$
The next step is to write $\vec{r}$ as a product of $r$ and of a
unit vector, $\vec{x}$, depending on angles. We put:
$$r_{\pm1}~=~r~x_{\pm1}\eqno(19)$$
$$r_0~=~r~x_0.\eqno(20)$$
It remains now to find a realization of $x_{\pm1}$ in terms of the
azimuthal angle $\varphi$ and of $x_0$, which is in fact equal to
$\cos\theta$, just as in $R(3)$ case. We found:
$$x_1~=~-e^{i\varphi}~\sqrt{q\over[2]}~\sqrt{1-q^2~x_0^2}~
q^{2N_0}\eqno(21)$$ 
$$x_{-1}~=~e^{-i\varphi}~\sqrt{1\over[2]q}~\sqrt{1-q^{-2}~x_0^2}~
q^{-2N_0}\eqno(22)$$ 
where the dilatation operator $N_0$ 
satisfying the relation
$$[~N_0~,~x_0^n~]~=~n~x_0^n\eqno(23)$$
and having the hermiticity property 
$$N_0^+~=~-N_0~-~1\eqno(24)$$
has been introduced in the expressions (21) and (22) in
order to fulfil the commutation relations (15) and (16).

Taking now into account the relations (19-24), and assuming
$$x_0^+~=~x_0\eqno(25)$$
we get for $x_\pm$ the normal hermiticity properties:
 $$x_1^+~=~{-1\over q}~x_{-1}\eqno(26)$$
$$x_{-1}^+~=~-q~x_1.\eqno(27)$$
All these arguments allow us to conclude that eqs.(19-23) define
the realization of 
the position vector $\vec{r}$ in $R_q(3)$ space.

The last step is to search for a realization of the $SU_q(2)$
generators. The expressions we propose for $L_+$ and $L_-$ are:
$$L_+~=~\sqrt{[2]}~e^{i\varphi}~{\tilde x}_1^{L_0+1}~{1\over
x_0}~{1-q^{-2N_0}\over 1-q^{-2}}~{\tilde x}_1^{-L_0}~q^{L_0}\eqno(28)$$
$$L_-~=~\sqrt{[2]}~e^{-i\varphi}~{\tilde x}_{-1}^{-L_0+1}~{1\over
x_0}~{1-q^{2N_0}\over 1-q^2}~{\tilde x}_{-1}^{L_0}~q^{L_0}\eqno(29)$$
where ${\tilde x}_{\pm1}={\rm e}^{\mp i\varphi}~x_{\pm1}$ depends on
$x_0$ only. Looking at the expressions (28) and (29) it becomes
clear why  the phase factor is removed from $x_{\pm1}$: expressions 
like $x_{\pm1}^{L_0}$ have no meaning,
while ${\tilde x}_{\pm1}^{L_0}$ is well defined.

Considering now the action of the operator $L_+$ on the non
normalized eigen
function of $L_0$ and of the Casimir operator C
$${\tilde Y}_{lm}(q,x_0,\phi)~=~{\rm e}^{im\varphi}~{\tilde x}_1^m
~\Theta_{lm}(x_0)\eqno(30)$$ 
we notice that ${\tilde x}_1^{-L_0}$ in (28) removes the factor
${\tilde x}_1^m$
in ${\tilde Y}_{lm}(q,x_0,\phi)$. In this way one prevents $q^{-2N_0}$
from acting on
$\tilde{x}_1$ and producing a troublesome result. The operator
$q^{-2N_0}$ in $L_+$ acts then on $\Theta_{lm}(x_0)$ only and 
e$^{i\varphi}~{\tilde x}_1^{L_0+1}$ creates the factor $x_1^m$
in right place.

In the well known $R(3)$ theory of angular momentum a different
mechanism prevents ${\partial\over\partial\theta}$ in $L_+$ from acting on
$x_1^m$: the term given by ${\partial\over\partial\theta}~x_1^m$ is exactly
cancelled by $i{\rm ctg}\theta{\partial\over\partial\varphi}~x_1^m$, 
finally remaining only the derivative ${\partial\over
\partial\theta}\Theta_{lm}(x_0)$.

It can be verified that the expressions (18), (28) and (29)
satisfy the commutation relations (1) and (2) and hence one can
conclude that they are the realization of the $SU_q(2)$
generators in $R_q(3)$ space. It can also be checked that the
position vector $\vec{r}$ defined in (19-22) behaves really like
a vector in this $SU_q(2)$ algebra, since it satisfies the
relations (6) and (7) with $L_\pm$ given by (28) and (29).
\vskip1cm
{\bf IV. THE LINEAR {\large$q$}-MOMENTUM {\large$\vec{p}$}}
\vskip1cm
In order to write down an expression for the linear momentum
$\vec{p}$, we separate it into a part perpendicular 
and another parallel to $\vec{x}$. The first one is defined with the
aid of the cross product $\vec{x}\times\vec{L}$ and the second
one is assumed to have the form $\vec{x}~{1\over
r}~f\left(r{d\over dr}+1\right)$, 
where $f$ is a function which will be defined in the following.
The components of the transverse part, denoted $\partial_k$, write
$$\partial_1~=~ q^{-1}~x_1~\Lambda_0~-~q~x_0~\Lambda_1~+~
x_1~c\eqno(31)$$
$$\partial_0~=~x_1~\Lambda_{-1}~-~\lambda~ x_0~\Lambda_0~-~x_{-1}
~\Lambda_1~+~x_0~c\eqno(32)$$
$$\partial_{-1}~=~- q~x_{-1}~\Lambda_0~+~q^{-1}~x_0~\Lambda_{-1}~+~
x_{-1}~c\eqno(33)$$
where $c$ is the invariant defined in eq.(12) and
the terms $x_k~c$ have been added to the cross product $\vec{x}\times
\vec{\Lambda}$ in order to ensure the well defined
character with respect to the hermitian conjugation operation
$$\partial^+_k~=~-\left(-{1\over q}\right)^k~\partial_{-k}. \eqno(34)$$
It can be checked that the quantities $\partial_k$ satisfy the
following relations: 
$$\partial_0\partial_{1}~=~q^{-2}~\partial_{1}\partial_0\eqno(35)$$
$$\partial_0\partial_{-1}~=~q^{2}~\partial_{-1}\partial_0\eqno(36)$$
$$\partial_1\partial_{-1}~=~\partial_{-1}\partial_1~+~\lambda^2~
\partial_0^2\eqno(37)$$
Eq.(35) has been obtained by commuting $\partial_0$ with
$\partial_1$, and eq.(36) is the hermitian conjugate of the
above one, while eq.(37) can be obtained either from eq.(35) or
(36) by using eq. (7).

Also, by multiplying equations (31-33) with the corresponding
$x_k$ and taking into account the commutation relations (15,16)
one gets:
$$\vec{x}~\vec{\partial}~=~-\vec{\partial}~\vec{x}=c~.\eqno(38)$$
By commuting the invariant c with $\vec{x}$ one finds:
$$\vec{\partial}~=~\lambda^{-2}~\left[c,\vec{x}\right].\eqno(39)$$
Taking now the matrix elements of the last relation one obtains:
$$\left\langle~
l+1~m'~\vert~\vec{\partial}~\vert~l~m\right\rangle~=~{[2l+2]\over[2l]}
~\left\langle~l+1~m'~\vert~\vec{x}~\vert~l~m\right\rangle\eqno(40)$$
$$\left\langle~
l-1~m'~\vert~\vec{\partial}~\vert~l~m\right\rangle~=~-{[2l]\over[2l]}
~\left\langle~l-1~m'~\vert~\vec{x}~\vert~l~m\right\rangle.\eqno(41)$$
From parity arguments one can also write:
$$\left\langle
l~m'\left\vert~\partial_k~\right\vert~l~m\right\rangle~=~0.\eqno(42)$$
By replacing the matrix elements of $\vec{\partial}$ with those
of $\vec{x}$ with the aid of eqs. (40) and (41) one can obtain
the eigenvalues of $\vec{\partial}^2$:
$$\left\langle~l~m\vert~\vec{\partial}^2~\vert~l~m\right\rangle~=~
-{[2l]\over[2]}~{[2l+1]\over[2]}~-~c_l^2.\eqno(43)$$
Taking into account all these relations, the realization we
found for the linear momentum $\vec{p}$ is:
$$\vec{p}~=~{-i\over r}~\left(~\vec{x}~(
r{\partial \over\partial r}~+~1)-\vec{\partial}~ \right).\eqno(44)$$
Then, in the $(2l+1)$ dimensional representation $\vec{p}~^2$
writes: 
$$\vec{p}~^2~=~-{1\over r}~{\partial\over\partial r}
\left(r{\partial\over\partial r}~+~1\right)~
+~{1\over r^2}\left({[2l]\over[2]}~{[2l+2]\over[2]}~+~c^2_l~-~c_l\right)
\eqno(45)$$
which in the limit $q=1$ becomes equal to the radial part of the
Laplace operator. 
\vskip1cm
{\bf V. THE {\large$q$}-ANGULAR MOMENTUM EIGEN FUNCTIONS}
\vskip0.5cm
The eigen vectors $\Phi_{lm}(q,x_0,\varphi)$ of the $(2l+1)$
dimensional irreducible representation of the $q$-angular
momentum are eigen functions of $L_0$ and of the Casimir
operator C. We begin by writing them as a polynomial in $x_0$
multiplied by $x_1^m$:
$$\Phi_{lm}(q,x_0,\varphi)~=~x_1^m~\sum_{k\geq0} ~a_k~x_0^k\eqno(46)$$
where the sum extends over odd and even $k$-values for
odd $(l-m)$ and respectively even $(l-m)$.

The equation:
$$L_+~L_-~\Phi_{lm}(q,x_0,\varphi)~=~[l+m]~
[l-m+1]~\Phi_{lm}(q,x_0,\varphi) \eqno(47)$$
gives the recursion relation:
$$a_{k+2}~=~-q^{-2m}~{[l-m-k]~[l+m+k+1]\over[k+1]~[k+2]}~a_k.\eqno(48)$$
For $(l-m)$ even we obtain: 
$$\Phi_{lm}(q,x_0,\varphi)~=~x_1^m~\left\{1-{[l-m][l+m+1]\over
[2]!}~\left(q^{-m}x_0\right)^2\right.$$
$$\left.+{[l-m][l-m-2][l+m+1][l+m+3]\over
[4]!}\left(q^{-m}x_0\right)^4-...\right\}\eqno(49)$$
while for $(l-m)$ odd we get:
$$\Phi_{lm}(q,x_0,\varphi)~=~x_1^m~\left\{{1\over[1]!}\left(
q^{-m}x_0\right)-{[l-m-1][l+m+2]\over
[3]!}~\left(q^{-m}x_0\right)^3\right.$$
$$\left.+{[l-m-1][l-m-3][l+m+2][l+m+4]\over
[5]!}\left(q^{-m}x_0\right)^5-...\right\}.\eqno(50)$$
In order to express these results in terms of a
$q$-hypergeometric series it is necessary to write all the
$q$-numbers $[n]$ in the form
$$[n]~=~{q^n-q^{-n}\over
q-q^{-1}}~=~[2]~{(q^2)^{n\over2}~-~(q^2)^{-{n\over2}}\over
q^2-q^{-2}}~=~ [2]~\left[{n\over2}\right]_{q^2}.\eqno(51)$$
For $(l-m)$ even we have then:
$$\Phi_{lm}(q,x_0,\varphi)~=~x_1^m~_2{\rm
F}_1\left(~q^2~;~{l+m+1\over2}~,~{-l+m\over2}~;~{1\over2}~;
~q^{-m}x_0^2~\right)\eqno(52)$$
while for $(l-m)$ odd we get: 
$$\Phi_{lm}(q,x_0,\varphi)~=~x_1^m~q^{-m}~x_0~_2{\rm
F}_1\left(~q^2~;~{l+m+2\over2}~,~{-l+m+1\over2}~;~{3\over2}~;
~q^{-m}x_0^2~\right).\eqno(53)$$
The index $q^2$ in $_2{\rm F}_1$ means that all the $q$-numbers
in the series development of $_2{\rm F}_1$ must be calculated
with $q^2$ instead of $q$.

We found that $\Phi_{lm}(q,x_0,\varphi)$ satisfies the following
simple relations:
$$x_1~{1\over
x_0}~{1-q^{-2N_0}\over1-q^{-2}}~\Phi_{lm}(q,x_0,\varphi)
~=~-[l-m]~[l+m+1]~\Phi_{l~m+1}(q,x_0,\varphi),\eqno(54)$$
for $(l-m)$ even, and
$$x_1~{1\over
x_0}~{1-q^{-2N_0}\over1-q^{-2}}~\Phi_{lm}(q,x_0,\varphi)
~=~\Phi_{l~m+1}(q,x_0,\varphi),\eqno(55)$$
for $(l-m)$ odd.

The normalized eigen functions are:
$$Y_{lm}(q,x_0,\varphi)=(-1)^{l-m\over2}\sqrt{[2l+1]\over4\pi}
\left({[l-m-1]!!\over[l-m]!!}~{[l+m-1]!!\over[l+m]!!}\right)^{1/2}
[2]^{m\over2}\Phi_{l~m}(q,x_0,\varphi),
\eqno(56)$$ 
for even $(l-m)$, and
$$Y_{lm}(q,x_0,\varphi)=(-1)^{l-m-1\over2}~\sqrt{[2l+1]\over4\pi}
\left({[l-m]!!\over[l-m-1]!!}{[l+m]!!\over[l+m-1]!!}\right)^{1/2}
[2]^{m\over2}\Phi_{l~m}(q,x_0,\varphi)\eqno(57)$$ 
for odd $(l-m)$ 
and the orthogonality relation is written as:
$$\int~Y_{l'm'}(q,x_0,\varphi)~Y_{lm}(q,x_0,\varphi)~d\varphi~
d[x_0]~=~ \delta_{ll'}~\delta_{mm'}\eqno(58)$$
where the integral over $\varphi$ is a normal one, while the
integrals over $d[x_0]$ are taken as:
$$\int^1_0~x_0^n~d[x_0]~=~{1\over[n+1]}.\eqno(59)$$
For the negative interval $(-1,0)$ from parity arguments we take:
$$\int_{-1}^0~x_0^n~d[x_0]~=~(-1)^n~{1\over[n]}.\eqno(60)$$
The relation (59) is in fact the result of a discrete
integration of $f(x_0)=x_0^n$, performed
by dividing the integration interval $(0,1)$ with an infinite
set of points $x_k=q^k$ for $q<1$ and writing 
$$\int_{-1}^0~f(x_0)~d[x_0]~=~
\sum_{k=0}^\infty~f(x_{2k+1})~(x_{2k}~-~x_{2k+2})~.\eqno(61)$$
Looking now for the properties of $Y_{lm}$, we found that, 
just as in the $R(3)$ case, the product $x_k~Y_{lm}$ can be
expressed in terms of $Y_{l\pm1~m+k}$ as follows:
$$x_1~Y_{lm}(q,x_0,\varphi)~=~q^{l-m}~\sqrt{[l+m+1][l+m+2]
\over[2][2l+1][2l+3]} 
~Y_{l+1~m+1}(q,x_0,\varphi)$$
$$-~q^{-l-m-1}~\sqrt{[l-m][l-m-1]\over[2][2l+1][2l-1]}~ 
Y_{l-1~m+1}(q,x_0,\varphi)\eqno(62)$$
\vskip0.5cm
$$x_0~Y_{lm}~=~q^{-m}~\sqrt{[l-m+1][l+m+1]\over[2l+1][2l+3]}
~Y_{l+1~m}(q,x_0,\varphi)$$
$$~-~q^{-m}~\sqrt{[l-m][l+m]\over[2l+1][2l-1]}~
Y_{l-1~m}(q,x_0,\varphi)\eqno(63)$$
\vskip0.5cm
$$x_{-1}~Y_{lm}(q,x_0,\varphi)~=~q^{l-m}~
\sqrt{[l-m+1][l-m+2]\over[2][2l+1][2l+3]} 
~Y_{l+1~m-1}(q,x_0,\varphi)$$
$$~-~q^{l-m+1}~\sqrt{[l+m][l+m-1]\over[2][2l+1][2l-1]}~ 
Y_{l-1~m-1}(q,x_0,\varphi).\eqno(64)$$
\vskip0.5cm
In addition, we have three relations which express the
noncommutativity of $x_k$ with $Y_{lm}$ and represent a
generalization of the equations ():
$$x_0~Y_{lm}(q,x_0,\varphi)~=~q^{-2m}~Y_{lm}(q,x_0,\varphi)~x_0
\eqno(65)$$ 
\vskip0.5cm
$$x_1~Y_{lm}(q,x_0,\varphi)~=~Y_{lm}(q,x_0,\varphi)~x_1$$
$$~+~{\lambda\over
\sqrt{[2]}}~q^{-m-1}~\sqrt{[l-m][l+m+1]}~Y_{l~m+1}(q,x_0,\varphi)~x_0
\eqno(66)$$
\vskip0.5cm
$$x_{-1}~Y_{lm}(q,x_0,\varphi)~=~Y_{lm}(q,x_0,\varphi)~x_{-1}$$
$$~-~{\lambda\over
\sqrt{[2]}}~q^{-m+1}~\sqrt{[l+m][l-m+1]}~Y_{l~m-1}(q,x_0,\varphi)~x_0.
\eqno(67)$$

The last two equations have been obtained from the first one
with the aid of $L_+$ and $L_-$ which rise and lower the index
$m$ of $Y_{lm}$. 
\vskip1cm
{\bf VI. {\large$q$}-DEFORMED SCHR\" ODINGER EQUATION}
\vskip0.5cm
Taking into account all the above results, we assume 
that the Hamiltonian entering the $q$- deformed Schrodinger
equation is: 
$${\cal H}~=~{1\over2}~\vec{p}~^2~+~V(r)\eqno(68)$$
where the operator $\vec{p}$ has been defined in the fourth
section. The eigen functions of this Hamiltonian are:
$$\Psi(r,x_0,\varphi)~=~r^L~u_L(r)~Y_{lm}(q,x_0,\varphi)\eqno(69)$$
where $L$ is the
solution of the  following equation:
$$L(L+1)~=~{[2l]\over[2]}{[2l+2]\over[2]}~+~c_l^2~-~c_l\eqno(70)$$
obtained from the condition that $u_L(r)$ is finite in the limit
$r\to0$.

This Schr\" odinger equation has simple solutions for Coulomb
potential $V(r)=-{1\over r}$ and for the oscillator potential
$V(r)={1\over2}~r^2$. The eigen values of the two Hamiltonians are:
$$E_{nl}~=~-{1\over2(n+L+1)^2}\eqno(71)$$
for the Coulomb potential and respectively:
$$E_{nl}~=~(2n+L+{3\over2}).\eqno(72)$$
for the oscillator potential,
$n$ being the radial quantum number and $L$ the solution of the equation
(70), usually not an integer. We notice that the spectrum is
degenerate with respect to the magnetic quantum number $m$.
The solution of the wave equation which does not depend on
$\theta$ and $\varphi$ gives for the mean value of $x_0^2$ the
value $R^2/[3]$ instead of $R^2/3$ obtained in the case of spherical
symmetry. It results then that the quadrupole momentum as well
as all the $2^{2n}$ poles are different from zero, although the
wave function does not depend on $\theta$ and $\varphi$.
This shows clearly that the Hamiltonian loosed the spherical
symmetry. One can mention however that it gained another, namely
the symmetry under the $SU_q(2)$ algebra.

Finally, we remark that there are three sources producing the
differences between the case of $q$-deformed Schr\" odinger
equation and the case of spherical symmetry: The first one is
that the $q$-harmonic wave function $Y_{lm}(q,x_0,\varphi)$
differs from the spherical one $Y_{lm}(\theta,\varphi)$. The
second reason is that 
the coefficient of the centrifugal potential in the radial
Schr\" odinger equation is $L(L+1)$, with $L$ given by eq.(70),
not $l(l+1)$, as in the sperical case. The third source is that
in the $q$-deformed case the integral over $x_0$ is performed 
according to the relations (59-60).
\vskip1cm
{\bf References}
\begin{enumerate}
\item Xing-Chang Song and Li Liao, J. Phys. A: Math. Gen. {\bf
25}, 623 (1992)
\item G. Rideau, P. Winternitz, J. Math. Phys. {\bf 34}, 6030, (1993)
\item Ya. I. Granovskii, A. S. Zhedanov, J. Phys. A: Math. Gen.
{\bf 26}, 4331 (1993)
\item Michele Irac-Astand Lett. Math. Phys. {\bf 36}, 169 (1996)
\item Denis Bonatsos, S. B. Drenska, P. P. Raychev, R. P.
Roussev and Yu. F. Smirnov, J. Phys. G: Nucl. Part. Phys. {\bf
17}, L67 (1991)
\item M. Nomura, J. Phys. Soc. Jpn. {\bf 59}, 439 (1990), V.
Rittenberg and M. Scheunert, J. Math. Phys. {\bf 33}, 436 (1992)

\end{enumerate}
\end{document}